\documentclass[12pt, reqno]{amsart}
\usepackage{xcolor}
\usepackage{hyperref}

\definecolor{DarkNavy}{HTML}{002060}    
\definecolor{Maroon}{HTML}{800000}      
\definecolor{PineGreen}{HTML}{004D20}   

\hypersetup{
    colorlinks=true,
    linkcolor=blue,   
    citecolor=blue,  
    urlcolor=blue       
}

\usepackage{amsmath,amsthm}
\usepackage{amsfonts}
\usepackage{amssymb}
\usepackage{comment}
\usepackage{graphicx}
\usepackage{xcolor}
\usepackage[left, pagewise, edtable]{lineno}
\usepackage{enumitem}
\usepackage[margin=1in,footskip=0.25in]{geometry}
\allowdisplaybreaks[4]   
\providecommand{\U}[1]{\protect\rule{.1in}{.1in}}

\theoremstyle{plain}
\newtheorem{theorem}{Theorem}[section] 

\newtheorem{corollary}[theorem]{Corollary}

\newtheorem{lemma}[theorem]{Lemma}

\theoremstyle{definition}

\theoremstyle{plain}

\title{On the Starshapedness of the Symmetrized Bidisc}
\author[Adith A.]{Adith A.}
\address{Adith A.\\Department of Physics, National Institute of Technology Calicut, 673601, Kozhikode, Keralam, India}
\email{adithammath@gmail.com}

\author[Geethika S.]{Geethika Sebastian}
\address{Geethika Sebastian\\Department of Mathematics, National Institute of Technology Calicut, 673601, Kozhikode, Keralam, India}
\email{geethikaseby@gmail.com}

\author[T. Banerjee]{Tathagata Banerjee}
\address{Tathagata Banerjee\\Department of Mathematics, Cochin University of Science and Technology, 682022, Cochin, Keralam, India}
\email{tathagatabanerjee85@gmail.com}

\date{}

\begin{document}

\begin{abstract}
    It is well known that the symmetrized bidisc $\mathbb{G}_2$ is star-like about the point $(0,0)$. In our attempt to check if there are other points about which $\mathbb{G}_2$ is star-like, we find that $\mathbb{G}_2$ cannot be starlike about any point $(s,p)\neq(0,0)\in \mathbb{G}_2$. 

\end{abstract}
\maketitle


\section{Introduction}\label{sec1}

The symmetrized bidisc and its geometric properties had already been a well-studied object as in \cite{MR2077158}, mainly for its importance from the function theoretic point of view as in the Nevalinna Pick interpolation theory. It has been shown in \cite{MR2985860} that the spectral version of the Nevalinna Pick interpolation problem has connections with the $n-$dimensional symmetrized polydisc. 
 
More importantly, the symmetrized bidisc is a non-convex domain that
exhibits behavior typically associated with convexity. For instance, it satisfies Lempert’s
Theorem  \cite{MR2077158}, a result previously known only for convex domains.
The geometric properties of the symmetrized bidisc, specially from the convexity point of view, are therefore important to study. Though the symmetrized bidisc is not a convex domain, it is star-like about the point $(0, 0)$ (\cite{MR2077158}Theorem 2.3). 

A set $M$ in $\mathbb{C}^n$ is called star-like/starshaped about the point $x_0\in{M}$ if for every point of ${M}$, the segment of the straight line from the point to $x_0$ lies in ${M}$.
In this article we prove the following result:
\begin{theorem}
   $\mathbb{G}_2$ is not star-like about any point $(s,p)\neq (0,0)\in \mathbb{G}_2 $.
    
\end{theorem}


We begin by fixing the notations to be used. Let $\mathbb{C}$ denote the space of all complex numbers. The unit disc in $\mathbb{C}$ is denoted by $$\mathbb{D}=\{z\in \mathbb{C}: |z|\leq 1\}.$$ The bidisc in $\mathbb{C}^2$ is denoted by $$\mathbb{D}^2=\{(z_1,z_2)\in \mathbb{C}^2: |z_1|<1, |z_2|<1\}.$$
The  symmetrized bidisc $\mathbb{G}_2$ is the image of the bidisc $\mathbb{D}^2$ under the ``symmetrization map'' 
$$\pi : (z_1,z_2) \mapsto (z_1 + z_2, \, z_1 z_2).$$ 
i.e., $$\mathbb{G}_2=\{(s,p)=(z_1+z_2,z_1z_2)\in \mathbb{C}^2:|z_1|<1, |z_2|<1)\}$$
The distinguised boundary in $\mathbb{C}^2$ is denoted by $$\Gamma(\mathbb{D}^2)=\{(e^{i\theta}, e^{i\phi}): \theta,\phi\in\mathbb{R}\}=\{(z_1,z_2)\in \mathbb{C}^2: |z_1|=1, |z_2|=1\}.$$ Consider the subset $T$ of $\Gamma(\mathbb{D}^2)$ given by $$T= \{(e^{i\theta}, e^{i\theta}): \theta\in\mathbb{R}\}.$$ The image of the set $T$ under $\pi$ will be denoted by the set
$$\pi(T)=\{(2e^{i\theta}, e^{2i\theta}): \theta\in\mathbb{R}\}$$ Note that $\pi(T)\subset \overline{\mathbb{G}}_2$, where $\overline{\mathbb{G}}_2=\{(z_1+z_2, z_1z_2)\in \mathbb{C}^2: |z_1|\leq 1, |z_2|\leq 1\}.$

For $x,y\in \mathbb{C}^2$, let $[x,y]$ denote the straight line joining the points $x$ and $y$. In particular $[x,y]=\{tx+(1-t)y, 0\leq t\leq 1\}.$

Now let us define the kernel of a set $M\subset \mathbb{C}^2$, denoted by
$$\text{ker}M=\{x\in M: [x,y]\subset M \;\text {for every}\; y\in M \}.$$ Theorem 1.2 in \cite{valentine1964convex}
states that the kernel of a set $M\subset \mathbb{C}^2$ is always convex.
\section{Proof of the main result}
We now proceed to the proof of the main result. The argument is built upon a sequence of intermediate results, which we establish first. Note that we are heavily going to use the following facts

\begin{itemize}
    \item $(s,p)\in \mathbb{G}_2$ iff $|s-\bar{s}p|<1-|p|^2$ (\cite{MR2077158} Theorem 2.1)
    \item $(s,p)\in \overline{\mathbb{G}_2}$ iff  $|s|\leq 2, |s-\bar{s}p|\leq 1-|p|^2$ (\cite{jarnicki2011invariant}, Lemma 7.1.4)  
\end{itemize}

\begin{lemma}\label{1}
    If $(s,p)\in \mathbb{G}_2$ then $(-s,p)\in \mathbb{G}_2.$
    \begin{proof}
        \begin{align*}
            (s,p)\in \mathbb{G}_2
            &\iff|s-\overline{s}p|<1-|p|^2\\
            &\iff|-(-s+\overline{s}p)|<1-|p|^2\\
            &\iff|-s-\overline{(-s)}p)|<1-|p|^2\\
            &\iff (-s,p)\in \mathbb{G}_2
        \end{align*}\qedhere
    \end{proof}
\end{lemma}

\begin{lemma}\label{3}
    If $(s,p)\in \ker\mathbb{G}_2$, then $(-s,p)\in \ker\mathbb{G}_2$.

    \begin{proof}
        We are given that $(s,p)\in \ker\mathbb{G}_2$, we need to prove that $(-s,p)\in \ker\mathbb{G}_2$.
        Suppose not, then there exists $(a,b)\in \mathbb{G}_2$ such that
$r(-s,p)+(1-r)(a,b)\notin \mathbb{G}_2$ for some $r\in(0,1).$
Thus for some $r\in(0,1)$, 
\begin{equation}\label{2}
 |(-rs+(1-r)a)-\overline{(-rs+(1-r)a)}(rp+(1-r)b)|\geq 1-|(rp+(1-r)b)|^2.   
\end{equation}

Now by Lemma \ref{1} we have that $(-a,b) \in \mathbb{G}_2.$
And since $(s,p)\in \ker\mathbb{G}_2$, thus
$r(s,p)+(1-r)(-a,b)\in \mathbb{G}_2$ for all $r\in(0,1).$ 
Therefore $|rs+(1-r)(-a)-\overline{(rs+(1-r)(-a))}(rp+(1-r)b)|< 1-|(rp+(1-r)b)|^2$. Multiplying the left hand side by $|-1|$ we get 
$$|-rs+(1-r)a-(-r\overline{s}+(1-r)\overline{a})(rp+(1-r)b)|<1-|rp+(1-r)b|^2.$$ This is a contradiction to (\ref{2}).\qedhere

    \end{proof}
\end{lemma}

The next corollary follows from  Lemma \ref{3} as $\ker\mathbb{G}_2$ is a convex set:
\begin{corollary}\label{11}
  If $(s,p)\in \ker\mathbb{G}_2$, then $(0,p)\in \ker\mathbb{G}_2$.
    
\end{corollary}

The proofs of the following two lemmas are analogous to those of Lemmas \ref{1} and \ref{3}; we include them here for completeness.

\begin{lemma}\label{4}
    If $(s,p)\in \mathbb{G}_2$ then $(is,-p)\in \mathbb{G}_2.$

    \begin{proof}
        \begin{align*}
            (s,p)\in \mathbb{G}_2\\
            &\iff|s-\overline{s}p|<1-|p|^2\\
&\iff|i(s-\overline{s}p)|<1-|p|^2\\
            &\iff|is-\overline{(is)}(-p)|<1-|p|^2\\
            &\iff (is,-p)\in \mathbb{G}_2
        \end{align*} \qedhere
    \end{proof}
\end{lemma}

\begin{lemma}\label{6}
        If $(s,p)\in \ker\mathbb{G}_2$, then $(-is,-p)\in \ker\mathbb{G}_2$.

        \begin{proof}
             We are given that $(s,p)\in \ker\mathbb{G}_2$, we need to prove that $(-is,-p)\in \ker\mathbb{G}_2$.
        Suppose not, then there exists $(a,b)\in \mathbb{G}_2$ such that
$r(-is,-p)+(1-r)(a,b)\notin \mathbb{G}_2$ for some $r\in(0,1).$
Therefore for some $r\in(0,1)$, 
\begin{equation}\label{5}
 |(-irs+(1-r)a)-\overline{(-irs+(1-r)a)}(-rp+(1-r)b)|\geq 1-|(-rp+(1-r)b)|^2.   
\end{equation}

Now by Lemma \ref{4} we have that $(ia,-b) \in \mathbb{G}_2.$
And since $(s,p)\in \ker\mathbb{G}_2$, thus
$r(s,p)+(1-r)(ia,-b)\in \mathbb{G}_2$ for all $r\in(0,1).$ 
i.e., $|rs+(1-r)(ia)-\overline{(rs+(1-r)(ia))}(rp+(1-r)b)|< 1-|(rp+(1-r)b)|^2$. Multiplying the left hand side of \ref{2} by $|-i|$ we get 
 $$|(-irs+(1-r)a)-\overline{(-irs+(1-r)a)}(-rp+(1-r)b)|< 1-|(-rp+(1-r)b)|^2.$$ This is a contradiction to (\ref{5}).\qedhere

        \end{proof}
\end{lemma}

The next corollary follows from  Lemma \ref{6} as $\ker\mathbb{G}_2$ is a convex set:

\begin{corollary}\label{9}
  If $(s,p)\in \ker\mathbb{G}_2$, then $(\frac{1-i}{2}s,0)\in \ker\mathbb{G}_2$.
    
\end{corollary}

\begin{theorem}\label{21}
   For $p\neq 0$ a point of the form $(0,p)\notin \ker\mathbb{G}_2$. 
   \begin{proof}
   Since $(0, p)\in \ker\mathbb{G}_2$ implies that $(0, p)\in \ker\overline{\mathbb{G}}_2$ by (\cite{valentine1964convex}, Theorem 1.10), we show by contradiction that $(0,p)\notin \ker\overline{\mathbb{G}}_2$.
      Let $p_1=(0,p)\in \ker\overline{\mathbb{G}}_2$
 and $p_2=(2e^{i\theta}, e^{2i\theta})\in \overline{\mathbb{G}}_2$. Let $p$ be such that $p=|p|e^{i\phi}$ and $\theta$ be such that $\theta\neq\frac{\phi}{2}$ and $\cos(\phi - 2\theta)> 0$. Consider $(a_t, b_t)=tp_2+(1-t)p_1=(2te^{i\theta}, te^{2i\theta}+(1-t)p)$, for $t\in [0, 1]$. We will show that if $(a_t,b_t)\in \overline{\mathbb{G}}_2 $, then there exists a $\delta<1$ such that $t\leq\delta.$  This is going to contradict the fact that if $(0,p)\in \ker \overline{\mathbb{G}}_2 $ then $(a_t,b_t)$ should belong to $\overline{\mathbb{G}}_2$ for every $t\in [0,1]$.  

Now, $(a_t, b_t)\in \overline{\mathbb{G}}_2$ implies $|a_t|\leq 2$ and $|a_t-\bar{a}_tb_t|\leq 1- |b_t|^2$ (\cite{jarnicki2011invariant}, Lemma 7.1.4). Clearly, $|a_t|=| 2te^{2i\theta}|\leq 2$ for all $t\in [0, 1]$. The second condition gives, $$|2te^{i\theta}-2te^{-i\theta}(te^{2i\theta}+(1-t)p)|\leq1-|te^{2i\theta}+(1-t)p|^2$$
  \begin{align*}
  &\implies |e^{-i\theta}||2te^{i\theta}-2te^{-i\theta}(te^{2i\theta}+(1-t)p)|\leq1-|te^{2i\theta}+(1-t)p|^2\\
   &\implies 2t|(1-t)-e^{-2i\theta}(1-t)p|\leq1-|te^{2i\theta}+(1-t)p|^2\\
   &\implies 2t(1-t)|1-pe^{-2i\theta}|\leq1-|te^{2i\theta}|^2-|(1-t)p|^2-2\Re (te^{-2i\theta}(1-t)p)\\
   &\implies 2t|1-pe^{-2i\theta}|\leq1+t-|p|^2+t|p|^2-2t\Re (e^{-2i\theta}p)\tag{1}
  \end{align*}
  Now consider $|1-pe^{-2i\theta}|^2=1+|p|^2-2\Re e^{-2i\theta}p$. Then the Binomial expansion for rational powers gives that

 \begin{equation}
  |1-pe^{-2i\theta}|=({1+|p|^2-2\Re (e^{-2i\theta}p
  )})^\frac{1}{2}=1+\frac{1}{2}(|p|^2-2\Re (e^{-2i\theta}p))+\epsilon \tag{2}
  \end{equation} where $\epsilon$ consists of the remainder terms in the Binomial expansion.\\
   Substituting (2) in (1) implies that $$t\leq\frac{1-|p|^2}{1+2\epsilon}=\delta\;\;  \text{say}$$
   Now if $\theta$ is chosen such that $\phi\neq 2\theta$ and $0<\cos(\phi-2\theta)$, then
   \begin{align*}
    &0<\cos(\phi-2\theta)<1, \\
    &\implies  |p|\cos(\phi-2\theta)<|p|\\
    &\implies \Re (pe^{-2i\theta})<|p|\;\;\;\;(\text{since}\; p=|p|e^{i\phi})\\
    &\implies (\Re (pe^{-2i\theta}))^2<|p|^2\\
     &\implies 1-2\Re (pe^{-2i\theta})+ (\Re (pe^{-2i\theta}))^2<1-2\Re (pe^{-2i\theta})+|p|^2\\
     &\implies (1-\Re (pe^{-2i\theta}))^2<1+\underbrace{|p|^2-2\Re (pe^{-2i\theta})}_x\\
     &\implies (1-\Re (pe^{-2i\theta}))<(1+x)^{\frac{1}{2}}\\
     &\implies \Re (pe^{-2i\theta})-1>-(1+x)^{\frac{1}{2}}\\
      &\implies |p|^2+ 2 \Re (pe^{-2i\theta})-2>-2(1+x)^{\frac{1}{2}}+|p|^2\\
      &\implies |p|^2>2+|p|^2-2 \Re (pe^{-2i\theta})-2(1+x)^{\frac{1}{2}}\\
      &\implies |p|^2>2+x-2(1+x)^{\frac{1}{2}}=2\{(1+\frac{1}{2}x)-(1+x)^\frac{1}{2}\}\\
      &\implies |p|^2>2|\epsilon|\footnotemark\\
      &\implies -|p|^2< 2\epsilon, \textrm{Since } \epsilon< 0\\
      &\implies 1 - |p|^2< 1 + 2\epsilon\\
      &\implies \delta=\frac{1-|p|^2}{1+2\epsilon}<1.
   \end{align*}
   \footnotetext{{Note that the Bernoulli's Inequality (\cite{bullen2003handbook}Section 2) states that: If $  x \geq -1, x \neq 0$, and if $0 < \alpha < 1$, then
$ (1+x)^{\alpha} < 1+\alpha x.$} Also $x=|p|^2-2\Re (pe^{-2i\theta})=|p|^2-2|p|\cos 2\theta>|p|^2+2|p|=|p|^2+2|p|+1-1=(|p|^2+1)^2-1>-1.$}
   \end{proof}

\end{theorem}
\begin{theorem}\label{22}
   For $s\neq 0$ a point of the form $(s,0)\notin \ker\mathbb{G}_2$. 
   \begin{proof}
    Let $s_1=(s,0)$
 and $s_2=(2e^{i\theta}, e^{2i\theta})$. Let $s$ be such that $s=|s|e^{i\phi}$ and $\theta$ be such that $\theta\neq\frac{\phi}{2}$. Consider $(a_t,b_t)=ts_2+(1-t)s_1=(2te^{i\theta}+(1-t)s, te^{2i\theta}).$ We will show that if $(a_t,b_t)\in \overline{\mathbb{G}}_2 $, then there exists a $\delta<1$ such that $t\leq\delta.$  This is going to contradict the fact that if $(s,0)\in \ker \mathbb{G}_2 $ then $(a_t,b_t)$ should belong to $\overline{\mathbb{G}}_2 $ for every $t\in [0,1]$.  

$(a_t, b_t)\in \overline{\mathbb{G}}_2$ implies $|a_t|\leq 2$ and $|a_t-\bar{a}_tb_t|\leq 1- |b_t|^2$ , i.e., $$|2te^{i\theta}+(1-t)s-\overline{(2te^{i\theta}+(1-t)s)}te^{2i\theta}|\leq1-|te^{2i\theta}|^2$$
  \begin{align*}
  &\implies |2te^{i\theta}+(1-t)s-(2te^{-i\theta}+(1-t)\bar{s})te^{2i\theta}|\leq1-t^2\\
  &\implies |2te^{i\theta}(1-t)+(1-t)(s-t\bar{s}e^{2i\theta})|\leq (1-t)(1+t)\\
  &\implies |2te^{i\theta}+s-t\bar{s}e^{2i\theta}|\leq 1+t\\
  &\implies |2te^{i\theta}+|s|e^{i\phi}-t|s|e^{-i\phi}e^{2i\theta}|\leq 1+t \tag{3}
  \end{align*}
Now consider $|2te^{i\theta}+|s|e^{i\phi}-t|s|e^{-i\phi}e^{2i\theta}|^2=|2te^{i\theta}+|s|(e^{i\phi}-te^{i(2\theta-\phi)})|^2$
\begin{align*}
    &=4t^2+|s|^2|e^{i\phi}-te^{i(2\theta-\phi)}| ^2+2\Re (2t|s|(e^{i(\phi-\theta)}-te^{i(\theta-\phi)}))
    \\
      &= 4t^2+|s|^2|1-te^{2i(\theta-\phi)}|^2+2\Re(2t|s|(e^{i(\phi-\theta)}-te^{i(\theta-\phi)}))\\
      &=4t^2+|s|^2|1-te^{2i\lambda}|^2+2\Re(2t|s|e^{-i\lambda}-2t^2|s|e^{i\lambda})\;\;\;\;(\text{where}\; \lambda =\theta-\phi)\\
      &= 4t^2+|s|^2(1+t^2-2t\cos 2\lambda)+4t|s|\cos \lambda-4t^2|s|\cos\lambda)\\
      \text {take }\lambda=\frac{\pi}{2}\\
      &= 4t^2+|s|^2(1+t^2+2t)\tag{4}
\end{align*}

Using (4) in (3) gives us $4t^2+|s|^2(1+t^2+2t)\leq (1+t)^2$,
\begin{align*}
&\implies  \frac{4t^2}{(1+t)^2}\leq1-|s|^2\\
&\implies \frac{2t}{1+t}\leq\sqrt{1-|s|^2}\\
&\implies \frac{1+t}{2t}\geq\frac{1}{\sqrt{1-|s|^2}}\\
&\implies \frac{1+t}{2t}-\frac{1}{2}\geq\frac{1}{\sqrt{1-|s|^2}}-\frac{1}{2}\\
&\implies 2t\leq\frac{1}{\frac{1}{\sqrt{1-|s|^2}}-\frac{1}{2}}\\
&\implies t\leq \frac{1}{2\left(\frac{1}{\sqrt{1-|s|^2}}-\frac{1}{2}\right)}=\delta \;(\text{say})
\end{align*}

Now since $s\neq 0$ \begin{align*}
    &\implies |s|^2>0\\
    &\implies 1-|s|^2<1
    \\
    &\implies\sqrt{1-|s|^2}<1\footnotemark\\
    &\implies \frac{1}{\sqrt{1-|s|^2}}>1\\
     &\implies \frac{2}{\sqrt{1-|s|^2}}>2\\
     &\implies \frac{2}{\sqrt{1-|s|^2}}-1>1\\
      &\implies \delta=\frac{1}{2\left(\frac{1}{\sqrt{1-|s|^2}}-\frac{1}{2}\right)} <1
     \end{align*}
     \footnotetext{Since $(s,0)\in\ker\mathbb{G}_2$ and since $\ker\mathbb{G}_2$ is convex we have $(rs,0)\in\ker\mathbb{G}_2$ for any $0<r<1$. Thus W.L.O.G., we can assume that $|s|<1.$ }
   \end{proof}

\end{theorem}
Hence, using Theorem \ref{21} and Theorem \ref{22}, we see that $(s, p)\notin \ker\mathbb{G}_2$ for $s\neq 0, p\neq 0$. Since otherwise, by Corollary \ref{11}, $(0, p)\in \ker \mathbb{G}_2$ which would contradict Theorem \ref{21}. Similarly, by Corollary \ref{9}, $(\frac{1-i}{2}s, 0)\in \ker\mathbb{G}_2$ which would contradict Theorem \ref{22}. Therefore we prove that the only point in the Kernel of $\mathbb{G}_2$ is $(0, 0)$.

\bibliographystyle{abbrv}
	\nocite{*}
\bibliography{sn-bibliography}
\end{document}